\documentclass[11pt,letterpaper]{article}

\usepackage[T1]{fontenc}
\usepackage{lmodern}
\usepackage[pdftex,margin=1in]{geometry}
\usepackage{enumerate}
\usepackage{amsmath}
\usepackage{amssymb}
\usepackage{fancyhdr}
\usepackage{lastpage}
\usepackage[alwaysadjust]{paralist}
\usepackage{amsthm}
\usepackage{graphicx}
\usepackage{wasysym}
\usepackage{xcolor}
\usepackage{listings}
\usepackage{hyperref}

\newtheorem{Thm}{Theorem}

\theoremstyle{definition}
\newtheorem{Def}[Thm]{Definition}

\theoremstyle{remark}

\theoremstyle{definition}

\theoremstyle{definition}

.

\DeclareRobustCommand{\maketitle}{%
	\begingroup
	\begin{center}
		\Large \textbf{Explicit Generating Functions for the Sum of the Areas Under
			Dyck and Motzkin Paths (and for Their Powers) }\\[\medskipamount]
		\large AJ Bu \\
	\end{center}
	\endgroup
}

\begin{document}
	\maketitle
	\begin{abstract}
		In this paper, we first describe how to find the generating function for the sum of the areas under generalized Dyck paths (with an arbitrary set of steps) using
		Motzkin paths as a motivating example. We then focus on Motzkin and Dyck paths, deriving  functional equations for them.
		We then  describe an algorithm to manipulate these functional equations for finding `perturbation expansions' of the solutions, with applications for deriving
		explicit generating functions for the sum of powers of areas under Dyck and Motzkin paths for any desired power.
	\end{abstract}
	
	\section{Introduction}
	
	Motzkin paths were extensively studied in  \cite{motzkin}, and the sum of powers of areas under Dyck paths were nicely covered by Robin Chapman \cite{RC}. However these
	papers and many others on these paths only  used human ingenuity. In this paper, we will demostrate the power of {\it symbolic computation} to get much further.
	
	We are interested in the generating functions for the sum of the areas under paths in the $xy-$plane, with a focus on Dyck and Motzkin paths.
	\begin{Def}
		A Motzkin path of length $n$ is a walk in the $xy-$plane from the origin $(0,0)$ to $(n,0)$ with atomic steps $U:=(1,1),$ $D:=(1,-1)$, and $F:=(1,0)$ that never goes below the $x-$axis.
	\end{Def}
	
	For example, the following are some Motzkin paths of length $4$. 
	\begin{align*}
		UDUD && UFFD && UFDF && FUDF && FFFF.
	\end{align*}
	The areas of these paths are $2,$ $4$, $3$, $1,$ and $0$, respectively.
	
	The bivariate weight enumerator for Motzkin paths with length $n$ and area $m$ satisfies the following functional equation
	$$M(x,q)=1+x M(x,q)+x^2 q M(qx,q)M(x,q).$$
	To prove this, let $\mathcal{M}$ denote the set of all Motzkin paths. Note that any path in $\mathcal{M}$ must fall into exactly one of the following cases -- the empty path, Motzkin paths that start with a flat step, or Motzkin paths that start with an up step. 
	
	If $M\in\mathcal{M}$ is the empty path, then it clearly has both area and length $0$. Thus the bivariate weight enumerator is $$m_0(q,x)=1$$ If $M$ begins with a flat step $F$, then we can write $$M=FM_0,$$
	where $M_0$ must also be a Motzkin path with the same area as $M$, since it still starts at height $0$. Thus, the bivariate weight enumerator for this case of Motzkin paths is  
	$$m_1(q,x)=xM(q,x)$$
	If $M$ begins with the step $U$, then let $D$ denote the first time $M$ returns to the $x-$axis and write
	\[M=U M_1 D M_0.\]
	$M_1$ must be a Motzkin path shifted to height $1$, and $M_0$ is a Motzkin path starting at height $0$. Since $M_0$ begins at height $0$, the area under the Motzkin path $M_0$ is the same as the area under the portion of $M$ it represents. Since $M_1$ is shifted to height 1, however, every step in $M_1$ has one more unit block below it. Thus, every step $x$ in $M_1$ must be replaced with $qx$ to get the correct area for that portion of $M$.  Since the extra $U$ and $D$ steps give a combined area of 1, the bivariate weight enumerator for Motzkin paths beginning with an up step is 
	\[m_2(q,x)=qx^2M(qx,q)M(q,x),\]
	resulting in the desired weight enumerator for all Motzkin paths.  
	\begin{Def}
		A Dyck path of length $n$ is a walk in the $xy-$plane from the origin $(0,0)$ to $(n,0)$ with atomic steps $U:=(1,1)$ and $D:=(1,-1)$ that never goes below the $x-$axis.
	\end{Def}
	Similarly, the bivariate weight enumerator for Dyck paths with length $n$ and area $m$ satisfies the following functional equation
	$$D(x,q)=1+x^2 q D(qx,q)D(x,q).$$
	\subsection{Maple Package}
	This article is accompanied by the Maple package \texttt{qEW.txt} and some sample outputs.  The accompanying files can be found at 
	\begin{center}
		\href{https://ajbu1.github.io/Papers/MotzArea/MotzArea.html}{https://ajbu1.github.io/Papers/MotzArea/MotzArea.html}.
	\end{center}
	
	\subsection{Using `Symbol-Crunching' to Find many terms of the area-weight-enumerators of generalized Dyck walks of length $n$}
	
	In addition to its independent interest, the {\it dynamical programming} approach to generate many terms of the sequence of area weight enumerators for {\it general} Dyck paths
	(with an {\it arbitrary} set of steps) will serve to confirm the correctness of the algorithm in the next section.
	
	The procedure $\texttt{qnwdK(S,K,q)}$ in our Maple package
	uses dynamical programming to find the enumerating function for the area of walks with steps in \texttt{S} of length $n=0,\dots,K$ that start and end at height $0$ and never have negative height.  
	Note that $S$ can be {\it any} set of steps, where the Dyck case is $S=\{[1,1],[1,-1]\}$ and the Motzkin case is  $S=\{[1,1],[1,-1], [1,0]\}$.
	
	First, consider $\mathcal{P}_{m,n}$, the set of paths of length $n\geq0$ with steps in $S$ that end at height $m\geq 0$ and never have negative height. Let $A_{m,n}(q)$ be the enumerating function for the area of paths in $\mathcal{P}_{m,n}$.
	
	Clearly, for $n=0$, the empty path gives an area of $0$.  Thus, the enumerating function is $$A_{m,0}(q)=1.$$ 
	For $n=1$, the only path that can end at height $m$ is the single step $\{(1,m)\}$, which has area $\frac{m}{2}$. Thus, 
	$$A_{m,1}(q)=\begin{cases}
		q^\frac{m}{2}, & (1,m)\in S\\
		0, & (1,m)\not\in S.
	\end{cases}$$ 
	For $n>1$, consider each possible final step for any path in $\mathcal{P}_{m,n}$.  A step $s\in S$ can be the last step if $m-s\geq 0$ and there exists a path $P$ of length $n-1$ with steps in $S$ that ends at height $m-s$ and never has a negative height. In other words,
	$$\mathcal{P}_{m,n}=\{Ps|s\in S,\text{ } m-s\geq 0,\text{ } P\in \mathcal{P}_{m-s,n-1}\}.$$
	The area under the last step $(1,s)$ is $\frac{2m-s}{2}$.  Thus, the weight enumerator for the area of paths in $\mathcal{P}_{m,n}$ is 
	$$A_{m,n}(q)=\sum_{\substack{s\in S \\m-s\geq 0}} q^\frac{2m-s}{2} A_{m-s,n-1}(q).$$
	This process is implemented in the procedure $\texttt{qnmwd(S,n,m,q)}$, which is then used in $\texttt{qnwdK(S,K,q)}$. 
	For example, looking at Motzkin paths,
	\[\texttt{qnwdK(\{[1,1],[1,0],[1,-1]\},5,q)}\]
	outputs 
	\[ \texttt{ $[1, 1, q + 1, q^2 + 2q + 1, q^4 + q^3 + 3q^2 + 3q + 1, q^6 + 2q^5 + 3q^4 + 4q^3 + 6q^2 + 4q + 1]$ .}\]
	Note that, to avoid negative height, any Motzkin path must end with $D=(1,-1)$ or $F=(1,0)$.  Since the paths end at height $0$, the area under these steps are $\frac{1}{2}$ and  $0$, respectively. Thus, $$A_{0,n}=q^\frac{1}{2}A_{1, n-1}+A_{0,n-1}.$$
	Breaking down the algorithm described above to find the first four terms of this outputted list,
	\begin{itemize}
		\item The only path of length $1$ is $[F]=[(1,0)]$, so $$A_{0,1}(q)=1.$$
		\item  Since $\mathcal{P}_{1,1}=\{[U]\}=\{[(1,0)]\}$, it follows that  $A_{1,1}(q)=q^\frac{1}{2}$. Thus,
		\begin{align*}
			A_{0,2}&=q^\frac{1}{2}A_{1,1}(q)+A_{0,1}(q)\\
			&=q+1.
		\end{align*}
		
		\item  For paths of length $3$ ending with $D=(1,-1)$, note that $$\mathcal{P}_{1,2}=\{FU,UF\}=\{[(1,0),(1,1)],[(1,1),(1,0)]\},$$ and so $A_{1,2}=q^\frac{1}{2}+q^\frac{3}{2}$. Thus,
		\begin{align*}
			A_{0,3}&=q^\frac{1}{2}A_{1,2}(q)+A_{0,2}(q)\\
			&=q^\frac{1}{2}(q^\frac{1}{2}+q^\frac{3}{2})+q+1\\
			&=q^2+2q+1.
		\end{align*}
		\item  For paths of length $4$, note that \begin{align*}\mathcal{W}_{1,3}&=\{FFU,FUF,UFF,UDU,UUD\}\\&=\{[(1,0),(1,0),(1,1)],[(1,0),(1,1),(1,0)],[(1,1),(1,0),(1,0)],\\
			&\hspace*{40pt}kL [(1,1),(1,-1),(1,1)],[(1,1),(1,1),(1,-1)]\}.   \end{align*} Therefore, $A_{1,3}=q^\frac{1}{2}+2q^\frac{3}{2}+q^\frac{5}{2}+q^\frac{7}{2}$, and
		\begin{align*}
			A_{0,4}&=q^\frac{1}{2}A_{1,3}(q)+A_{0,3}(q)\\
			&=q^4+q^3+3q^2+3q+1.
		\end{align*}
		
	\end{itemize}

	\section{Finding $\frac{d^k}{dq^k}\left[f(x,q) \right]\big|_{q=1}$}
	If the weight enumerator of a set of paths is satisfied by the following functional equation
	$$f(x,q)=P(x,q)+Q(x,q)f(x,q)+R(x,q)f(x,q)f(qx,q)$$
	for some given bivariate polynomials $P(x,q)$, $Q(x,q)$, and $R(x,q),$ then plugging in $q=1$ gives
	$$f(x,1)=P(x,1)+Q(x,1)f(x,1)+R(x,1)f(x,1)^2,$$
	which we can use to solve for $f(x,1).$  The order $n$ Taylor polynomial of $f(x,q)$ about $q=1$ satisfies
	$$\sum_{k=0}^n \frac{(q-1)^k}{k!}f^{(k)}(x,1)=P+Q \sum_{k=0}^n \frac{(q-1)^k}{k!}f^{(k)}(x,1)+ R \sum_{k=0}^n \frac{(q-1)^k}{k!}f^{(k)}(x,1) \sum_{k=0}^n \frac{(q-1)^k}{k!}f^{(k)}(qx,1),$$
	where $f^{(k)}(x,q)=\frac{d^k}{dq^k} f(x,q).$
	Looking at the coefficient of $(q-1)^k$, we can express  $f^{(k)}(x,1)$ as the sum of derivatives $f^{(\ell)}(x,1)$ where $\ell<k$ and derivatives of functions of $x$ with respect to $x$.  Since we have an expression for $f(x,1)$, we can simply compute any order derivative with respect to $x$ as well as $f_q(x,1)$.  Thus, to find $f^{(n)}(x,1)$, we can repeat this process with the coefficient of $f^{(k)}(x,1)$ for $k=1,\dots, n.$
	
	This process is implemented by the procedure \texttt{DerK(P,Q,R,q,x,K,f)}, which outputs a list whose $k$-th entry is $\frac{d^{k-1}}{dq^{k-1}}\left[f(x,q) \right]\big|_{q=1}$. Rather than outputting algebraic equations as seen in \cite{bz}, this procedure produces closed-form expressions in terms of radicals.

	\subsection{Motzkin Paths}
	As previously noted, the Motzkin paths satisfy the following functional equation
	$$M(x,q)=1+x M(x,q)+x^2 q M(qx,q)M(x,q).$$
	Solving this functional equation for $q=1$, we get that
	\begin{align*}
		M(x,1)&=\frac{1-x+\sqrt{-3x^2-2x+1}}{2x^2} &\text{or}&& M(x,1)&=\frac{1-x-\sqrt{-3x^2-2x+1}}{2x^2}.
	\end{align*}
	Since only the second equation has a Taylor series expansion about $x=0$, we know that this is $M(x,1)$.  Now, for finding the first derivative, note that 
	$$\sum_{k=0}^n \frac{(q-1)^k}{k!} M^{(k)}(x,1)=1+x\sum_{k=0}^n \frac{(q-1)^k}{k!} M^{(k)}(x,1) + q x^2 \sum_{k=0}^n \frac{(q-1)^k}{k!} M^{(k)}(x,1) \sum_{k=0}^n \frac{(q-1)^k}{k!} M^{(k)}(qx,1)$$
	The coefficient of $q-1$ on both sides give  us 
	$$M_q(x,1)=xM_q(x,1) + x^2M(x,1)\bigg(xM_x(x,1)+2M_q(x,1) +M(x,1)\bigg).$$
	Therefore,
	$$M_q(x,1)=\frac{x^3M(x,1)M_x(x,1)+x^2 M^2 (x,1) }{1-x-2x^2 M(x,1)}.$$
	Plugging in $M(x,1)=\frac{1-x-\sqrt{-3x^2-2x+1}}{2x^2},$ we get
	$$M_q(x,1)=\frac{\left(x-1+\sqrt{-3x^2-2x+1}\right)^2}{4x^2(-3x^2-2x-1)}$$
	To find $M^{(n)}(x,1)$, we can repeat this process with the coefficient of $M^{(k)}(x,1)$ for $k\leq n.$
	
	In a little over 2 seconds,
	\[\texttt{DerK(1,x,x$^2$*q,q,x,10,f)},\]
	can output the list whose entries are $M^{(k)}(q,1):=\frac{d^k}{dq^k}\left[M(x,q)\right]|_{q=1}$ for $k=0,\dots,10$. For example, looking at the first two terms of the output, we have
	\begin{align*}
		M(x,1)=&\frac{1-x-\sqrt{-3x^2-2x+1}}{2x^2}, &\text{and}\\
		M_q(x,1)=&\frac{\left( 1- x - \sqrt{-3x^2 - 2x + 1}\right)^2}{4x^2(-3x^2-2x+1)}
	\end{align*}
	The Maclaurin Series of $M(x,1)$ is
	$$1 + x + 2x^2 + 4x^3 + 9x^4 + 21x^5 + 51x^6 + 127x^7 + 323x^8 + 835x^9 + 2188x^{10} + 5798x^{11} + 15511x^{12} + O(x^{13}),$$
	and it is the weight enumerator of the number of Motzkin paths of length $n$, which is A001006 on \cite{oeis}, \href{https://oeis.org/A001006}{https://oeis.org/A001006}.  
	The Maclaurin series of $M_q(x,1)$  is
	$$x^2 + 4x^3 + 16x^4 + 56x^5 + 190x^6 + 624x^7 + 2014x^8 + 6412x^9 + 20219x^{10} + 63284x^{11} + 196938x^{12} + O(x^{13}),$$ 
	which is the weight enumerator of the total area under all Motzkin paths of length $n$ and A057585 on \cite{oeis}, \href{https://oeis.org/A057585}{https://oeis.org/A057585}.
	
	We also get higher factorial moments. For example, \\
	$M_{qq}(x,1)=1/2(6(-3x^2-2x+1)^{1/2}x^2+9x^2-(-3x^2-2x+1)^{1/2}x+6x+3(-3x^2-2x+
	1)^{1/2}-3)(-1+x+(-3x^2-2x+1)^{1/2})/(3x^2+2x-1)^3,$
	\begin{center} and \end{center}
	$M_{qqq}(x,1)=-3/2(9(-3
	x^2-2x+1)^{1/2}x^4-9x^5+18(-3x^2-2x+1)^{1/2}x^3+51x^4-23(-3x^2-2x+1
	)^{1/2}x^2-19x^3+4(-3x^2-2x+1)^{1/2}x+29x^2-4(-3x^2-2x+1)^{1/2}-8x+4
	)(x-1+(-3x^2-2x+1)^{1/2})/(3x^2+2x-1)^4.$
	
	The Maclaurin series of $M_{qq}(x,1)$ is
	$$2x^3 + 24x^4 + 142x^5 + 720x^6 + 3224x^7 + 13478x^8 + 53508x^9 + 204698x^{10} + O(x^{11}),$$ and the  Maclaurin series of $M_{qqq}(x,1)$ is
	$$30x^4 + 336x^5 + 2742x^6 + 17268x^7 + 95388x^8 + 477900x^9 + 2235876x^{10}  + O(x^{11}).$$
	The weight enumerator for the sum of the squares of the areas of Motzkin paths of length $n$ is given by the Maclaurin series of $M_{qq}(x,1)+M_q(x,1),$
	$$x^2 + 6x^3 + 40x^4 + 198x^5 + 910x^6 + 3848x^7 + 15492x^8 + 59920x^9 + 224917x^{10} + O(x^{11}),$$ 
	which is A367778 on \cite{oeis}, \href{https://oeis.org/A367778}{https://oeis.org/A367778}.  
	
	The weight enumerator for the sum of the cubes of the areas of Motzkin paths of length $n$ is given by the Maclaurin series of $M_{qqq}(x,1)+3M_{qq}(x,1)+M_q(x,1),$
	$$x^2 + 10x^3 + 118x^4 + 818x^5 + 5092x^6 + 27564x^7 + 137836x^8 + 644836x^9 + 2870189x^{10} + O(x^{11}),$$
	which is A367779 on \cite{oeis}, \href{https://oeis.org/A367779}{https://oeis.org/A367779}.
	
	\subsection{Dyck Paths}
	Looking at Dyck paths, we input 
	$$\texttt{DerK(1,0,x$^2$*q,q,x,10,f)}.$$
	The first four terms of the output gives
	\begin{align*}
		D(x,1)&=\frac{1-\sqrt{1-4x^2}}{2x^2}\\
		D_q(x,1)&=\frac{(1-\sqrt{1-4x^2})^2}{16x^4-4x^2}\\
		D_{qq}(x,1)&=\frac{(8x^2\sqrt{1-4x^2}+12x^2+3\sqrt{1-4x^2}-3)(-1+\sqrt{1-4x^2})}{2}\\
		D_{qqq}(x,1)&=\frac{-6(4x^4\sqrt{-4x^2+1}+16x^4-7x^2\sqrt{-4x^2+1}+7x^2-\sqrt{-4x^2+1}
			+1)(-1+\sqrt{-4x^2+1})}{(4x^2-1)^4}
	\end{align*}
	The Maclaurin series of $D(x,1)$ is 
	$$1 + x^2 + 2x^4 + 5x^6 + 14x^8 + 42x^{10} + 132x^{12} + 429x^{14} + 1430x^{16} + O(x^{18}),$$ which is the weight enumerator of all Dyck paths of length $n$ and A000108 on \cite{oeis}, \href{https://oeis.org/A000108}{https://oeis.org/A000108}. 
	The Maclaurin series of $D_q(x,1)$ is 
	$$x^2 + 6x^4 + 29x^6 + 130x^8 + 562x^{10} + 2380x^{12} + 9949x^{14} + 41226x^{16} + O(x^{18}),$$ the weight enumerator for the total area of all Dyck paths of length $n$, which is A008549 on \cite{oeis}, \href{https://oeis.org/A008549}{https://oeis.org/A008549}. \\
	The Maclaurin series of $D_{qq}(x,1)$ is
	$$14x^4 + 160x^6 + 1226x^8 + 7864x^{10} + 45564x^{12} + 247136x^{14} + 1279810x^{16} + 6404424x^{18} + O(x^{20}),$$
	and the Maclaurin series of $D_{qqq}(x,1)$ is
	$$24x^4 + 840x^6 + 11736x^8 + 114744x^{10} + 922224x^{12} + 6541776x^{14} + 42543480x^{16} + 259525464x^{18} + O(x^{20}).$$
	The weight enumerator for the sum of the squares of the areas of Dyck paths of length $n$ is given by the Maclaurin series of $D_{qq}(x,1)+D_q(x,1),$
	$$x^2 + 20x^4 + 189x^6 + 1356x^8 + 8426x^{10} + 47944x^{12} + 257085x^{14} + 1321036x^{16}+ O(x^{18}),$$
	which is A367780 on \cite{oeis}, \href{https://oeis.org/A367780}{https://oeis.org/A367780}.

	The weight enumerator for the sum of the cubes of the areas of Dyck paths of length $n$ is given by the Maclaurin series of $D_{qqq}(x,1)+3D_{qq}(x,1)+D_q(x,1),$
	$$x^2 + 72x^4 + 1349x^6 + 15544x^8 + 138898x^{10} + 1061296x^{12} + 7293133x^{14} + 46424136x^{16} + O(x^{18}).$$
	This does not appear on OEIS as of September 12, 2023.
	\section{Conclusion}
	In this paper, we demonstrate how to use dynamical programming to find the weight enumerator for the area paths of length $n$ with steps in a given set $S$ that start and end at height $0$ and never have negative height.  We also describe how to find the weight enumerator for such paths when, instead of a set of steps $S$, we are given bivariate polynomials $P(x,q)$, $Q(x,q)$, and $R(x,q)$ such that the weight enumerator $f(x,q)$ satisfies $$f(x,q)=P(x,q)+Q(x,q)f(x,q)+R(x,q)f(x,q)f(xq,q).$$ We then present a method for finding $\displaystyle f^{(k)}(x,1):=\frac{d^k}{dq^k}[f(x,q)]\big|_{q=1}$.
	
	These methods are fully automated in the accompanying Maple package \texttt{qEW.txt}, displaying how the power of computer algebra and using calculus allows us to generate quite a few moments. In the paper, we demonstrate these methods with the bivariate weight enumerators for both Motzkin paths and Dyck paths with length $n$ and area $m$.  Moreover, we show how these procedures can be used to produce the Maclaurin series of $\frac{d^k}{dq^k}\left[f(x,q) \right]\big|_{q=1}$, allowing us to find the generating function for the total area under such paths of length $n$  as well as for the sum of a given power of the areas.
	
	For further study, we can look at the average areas and the variance.  Given a family of paths, let $a_0(n)$ be the number of such paths of length $n$, $a_1(n)$ be the total area under such paths of length $n$, and $a_2(n)$ be the sum of the squares of the areas under such paths of length $n$.  Using the accompanying Maple package \texttt{qEW.txt}, we can generate 10,000 (or more) terms of the sequences of the average areas $\displaystyle \left\{\frac{a_1(n)}{a_0(n)}\right\}$ and the variances $\displaystyle \left\{\frac{a_2(n)}{a_0(n)} -\left(\frac{a_1(n)}{a_0(n)}\right)^2\right\}$ and use numerics for the asymptotics.

	\vspace{11pt}
	\subsection*{Acknowledgements}
	Thank you to Doron Zeilberger for helpful feedback and guidance in research for this paper.

	%\makefoot
\end{document}